\documentclass[graybox]{svmult}

\usepackage{mathrsfs}
\usepackage{mathtools}
\usepackage{todonotes}

\usepackage{pgf,tikz}
\usepackage{mathrsfs}
\usetikzlibrary{arrows}

\usepackage{mathptmx}       
\usepackage{helvet}         
\usepackage{courier}        
\usepackage{type1cm}        
\usepackage{makeidx}         
\usepackage{graphicx}        
\usepackage{multicol}        
\usepackage[bottom]{footmisc}

\usepackage{cite}

\usepackage{amsmath,amssymb,bm}

\makeindex             

\begin{document}
 
\title*{High-Order Isogeometric Methods for Compressible Flows.}
\subtitle{II. Compressible Euler Equations}
\titlerunning{High-Order Isogeometric Methods for Compressible Flows II.}
\author{Matthias M{\"o}ller and Andrzej Jaeschke}
\institute{
Matthias M{\"o}ller \at Delft University of Technology, Faculty of Electrical Engineering, Mathematics and Computer Science, Delft Institute of Applied Mathematics, Van Mourik Broekmanweg 6, 2628 XE Delft, The Netherlands, \email{m.moller@tudelft.nl}
\and
Andrzej Jaeschke \at \L{}\'od\'z University of Technology, Institute of Turbomachinery, ul. W\'olcza\'nska 219/223, 90-924 \L{}\'od\'z, Poland, \email{andrzej.jaeschke@p.lodz.pl}
}
%

%
\maketitle

\abstract*{
This work extends the high-resolution isogeometric analysis approach established in \cite{igaAFC1_fef2017} to the equations of gas dynamics. The group finite element formulation is adopted to obtain an efficient assembly procedure for the standard Galerkin approximation, which is stabilized by adding artificial viscosities proportional to the spectral radius of the Roe-averaged flux-Jacobian matrix. Excess stabilization is removed in regions with smooth flow profiles with the aid of algebraic flux correction \cite{KBNII}. The underlying principles are reviewed and it is shown that FCT-type flux limiting originally derived for nodal low-order finite elements ensures positivity-preservation for high-order B-Spline discretizations.
}

\abstract{
This work extends the high-resolution isogeometric analysis approach established in \cite{igaAFC1_fef2017} to the equations of gas dynamics. The group finite element formulation is adopted to obtain an efficient assembly procedure for the standard Galerkin approximation, which is stabilized by adding artificial viscosities proportional to the spectral radius of the Roe-averaged flux-Jacobian matrix. Excess stabilization is removed in regions with smooth flow profiles with the aid of algebraic flux correction \cite{KBNII}. The underlying principles are reviewed and it is shown that linearized FCT-type flux limiting \cite{Kuzmin2009} originally derived for nodal low-order finite elements ensures positivity-preservation for high-order B-Spline discretizations.
}

\section{Introduction}
\label{sec:introduction}

Compressible fluid flow problems have traditionally been solved by low-order finite element and finite volume schemes, which were equipped with stabilization techniques like, e.g., SUPG, FCT/TVD and nonlinear shock-capturing, in order to resolve flow patterns like shock waves and contact discontinuities without producing nonphysical undershoots and overshoots. A recent trend in industrial CFD applications especially involving turbulent flows is the use of high-order methods, which enable a more accurate representation of curved geometries and thus a better resolution of boundary layers and, in general, provide better accuracy per degree of freedom (DOF). Moreover, high-order methods have a favorable compute-to-data ratio, which makes them particularly attractive for use on modern high-performance computing platforms, where memory transfers are typically the main bottleneck.

Despite the huge success of isogeometric analysis (IGA) \cite{IGABook} in structural mechanics, incompressible fluid mechanics and fluid-structure interaction, publications on its successful application to \emph{compressible} flows are rare \cite{Jaeschke,DUVIGNEAU2018443,IGAart2}. This may be attributed to the challenges encountered in developing shock-capturing techniques for \emph{continuous} high-order finite element methods. In this work we present an isogeometric approach for solving the compressible Euler equations within an IGA framework, thereby adopting the concept of algebraic flux correction (AFC) as stabilization technique \cite{KBNI,KBNII}. In particular, we show that the linearized FCT-type limiting strategy introduced in \cite{Kuzmin2009,Kuzmin2010} carries over to B-Spline discretizations.

\section{High-resolution isogeometric analysis}

This section briefly describes the design principles of our IGA-AFC scheme and highlights the novelties and main differences to its nodal finite element counterpart. For a comprehensive description of FEM-AFC the reader is referred to \cite{KBNI,KBNII}.

\subsection{Governing equations}
\label{sec:goveq}

Consider the $d$-dimensional compressible Euler equations in divergence form
\begin{equation}
\frac{ \partial U}{\partial t} + \nabla \cdot \mathbf{F}(U) = 0.
\label{eq:conslaw}
\end{equation}
Here, $U:\mathbb{R}^{d}\to\mathbb{R}^{d+2}$ denotes the state vector of conservative variables
\begin{equation}
U 
= \begin{bmatrix}
U_1,\dots,U_{d+2}
\end{bmatrix}^\top
= \begin{bmatrix}
\rho,
\rho \boldsymbol{\mathrm{v}},
\rho E
\end{bmatrix}^\top,
\label{eq:statevector}
\end{equation}
and $\mathbf{F}:\mathbb{R}^{d+2}\to \mathbb{R}^{(d+2)\times d}$ stands for the tensor of inviscid fluxes
\begin{equation}
\mathbf{F} 
= \begin{bmatrix}
F^1_1 & \hdots & F^d_1\\
\vdots & \ddots & \vdots\\
F^1_{d+2} & \hdots & F^d_{d+2}
\end{bmatrix}
= \begin{bmatrix}
\rho \boldsymbol{\mathrm{v}}\\
\rho \boldsymbol{\mathrm{v}}\otimes \boldsymbol{\mathrm{v}}+  p I\\
\rho E \boldsymbol{\mathrm{v}} + p \boldsymbol{\mathrm{v}}
\end{bmatrix}
\label{eq:inviscidflux}
\end{equation}
with density $\rho$, velocity $\mathbf{v}$, total energy $E$, and $I$ denoting the $d$-dimensional identity tensor. For an ideal polytropic gas, the pressure $p$ is given by the equation of state
\begin{equation}
p = (\gamma -1) \left(\rho E - 0.5\rho|\boldsymbol{\mathrm{v}}|^2\right),
\label{eq:eos}
\end{equation}
 where $\gamma$ denotes the heat capacity ratio, which equals $\gamma = 1.4$ for dry air.
The governing equations are equipped with initial conditions prescribed at time $t=0$ 
\begin{equation}
U(\boldsymbol{\mathrm{x}}, 0) = U_0 (\boldsymbol{\mathrm{x}}) \;\; \text{in } \Omega,
\label{eq:initcond}
\end{equation}
and boundary conditions of Dirichlet and Neumann type, respectively
\begin{equation}
U = G(U,U_\infty) \;\; \text{on} \; \varGamma_D,
\quad
\mathbf{n} \cdot \mathbf{F} = F_n(U,U_\infty) \;\; \text{on}\; \varGamma_N.
\label{eq:bdrcond}
\end{equation}
Here, $\mathbf{n}$ is the outward unit normal vector and $U_\infty$ denotes the vector of 'free stream' solution values, which are calculated and imposed as outlined in \cite{KBNII}.

\subsection{Spatial discretization by isogeometric analysis}
\label{sec:spatial}

Application of the Galerkin method to the variational form of the first-order conservation law system \eqref{eq:conslaw} yields the following system of semi-discrete equations \cite{KBNII}
\begin{equation}
\sum_{j} \left(\int_{\Omega}\varphi_i \varphi_j \mathrm{d} \mathbf{x} \right ) \frac{d U_j}{d t} -
\sum_{j} \left(\int_{\Omega}\nabla\varphi_i \varphi_j \mathrm{d} \mathbf{x} \right ) \cdot \mathbf{F}_j +
\int_{\Gamma_n}\varphi_i F_n\,\mathrm{d}s= 0,
\label{eq:semidisc}
\end{equation}
where $\mathbf{F}_j(t):=\mathbf{F}(U_j(t))$ denotes the value of \eqref{eq:inviscidflux} evaluated at the $j$-th solution coefficient $U_j(t)$ at time $t$. This approach is known as Fletcher's group formulation \cite{Fletcher}, which amounts to expanding $U^h\approx U$ and $\mathbf{F}^h\approx\mathbf{F}$ into the same basis
$\{\varphi_j\}$, that is
\begin{equation}
U^h(\mathbf{x},t)=\sum_{j}U_j(t)\varphi_j(\mathbf{x}),
\quad
\mathbf{F}^h(\mathbf{x},t)=\sum_{j}\mathbf{F}_j(t)\varphi_j(\mathbf{x}).
\label{eq:fletcher}
\end{equation}
To further simplify the notation, let us define the consistent mass matrix $M_C:=\{m_{ij}\}$ and the discretized divergence operator $\mathbf{C}:=\{\mathbf{c}_{ij}\}$ as follows
\begin{equation}
m_{ij}=\int_\Omega\varphi_i\varphi_j\,\mathrm{d}\mathbf{x},\qquad
\mathbf{c}_{ij}=\int_\Omega\varphi_i\nabla\varphi_j\,\mathrm{d}\mathbf{x}.
\label{eq:coeffmat}
\end{equation}
Then the semi-discrete system \eqref{eq:semidisc} can be written in compact matrix form as
\begin{equation}
M_C\frac{dU_k}{dt}-\sum_{l=1}^d \left[C^l\right]^\top F^l_k+S_k(U)=0,
\label{eq:semidisc_spmv}
\end{equation}
where superscript $l=1,\dots,d$ refers to the $l$-th spatial component of the discrete divergence operator $\mathbf{C}$ and the tensor of inviscid fluxes $\mathbf{F}$, respectively, and subscript $k=1,\dots,d+2$ stands for the component that corresponds to the $k$-th variable. Here, $S_k(U)$ accounts for the contribution of boundary fluxes; see \cite{KBNII} for more details.

As explained in Sections 2.1 and 2.2 of \cite{igaAFC1_fef2017}, we adopt tensor-product B-Spline basis functions for approximating the numerical solution \eqref{eq:fletcher} and modeling the  domain $\Omega$. As a consequence, the integrals in \eqref{eq:coeffmat} are evaluated by introducing the 'pull-back' operator $\phi^{-1}:\Omega\to\hat\Omega$ and applying numerical quadrature on the reference domain $\hat\Omega=[0,1]^d$ as it is common practice in the IGA community.

%

\subsection{Temporal discretization by explicit Runge-Kutta methods}
\label{sec:time}

The semi-discrete system \eqref{eq:semidisc_spmv} is discretized in time by an explicit strong stability preserving (SSP) Runge-Kutta time integration schemes  of order 
three \cite{Gottlieb}
\begin{align}
MU^{(1)}&=MU^n+\Delta t R(U^n)\label{eq:rk3_1}\\
MU^{(2)}&=\frac34 MU^n+\frac14\left(MU^{(1)}+\Delta t R(U^{(1)})\right)\label{eq:rk3_2}\\
MU^{n+1}&=\frac13 MU^n+\frac23\left(MU^{(2)}+\Delta t R(U^{(2)})\right),\label{eq:rk3_3}
\end{align}
where $M=I\otimes M_C$ is the block counterpart of the mass matrix and $R(U)$ represents all remaining terms involving inviscid and boundary fluxes. If $M_C$ is replaced by its row-sum lumped counterpart (see below) then the above Runge-Kutta schemes reduce to scaling the right-hand sides by the inverse of a diagonal matrix.

\subsection{Algebraic flux correction}
\label{sec:afc}

The Galerkin method \eqref{eq:semidisc} is turned into a high-resolution scheme by applying algebraic flux correction of linearized FCT-type \cite{Kuzmin2009}, thereby adopting the primitive variable limiter introduced in \cite{Kuzmin2010}. The description is intentionally kept short and addresses mainly the extensions of the core components of FEM-AFC to IGA-AFC.

\textbf{Row-sum mass lumping}. A key ingredient in all flux correction algorithms for time-dependent problems is the decoupling of the unknowns in the transient term of \eqref{eq:semidisc} by performing row-sum mass lumping, which turned out to be one of the main problems in generalizing FEM-AFC to higher order nodal Lagrange finite elements since the presence of negative off-diagonal entries leads to singular matrices. The following theorem shows that IGA-AFC is free of this problem by design.

\begin{theorem}
All diagonal entries and the non-zero off-diagonal entries of the consistent mass matrix $M_C$ are strictly positive if tensor-product B-Spline basis functions ${\varphi_j}$ are adopted. Hence, the row-sum lumped mass matrix $M_L:=\text{diag}(m_i)$ is unconditionally invertible with strictly positive diagonal entries
\begin{equation}
m_i:=
\sum_{j} \int_{\Omega}\varphi_i(\mathbf{x})\varphi_j(\mathbf{x})\,\mathrm{d}\mathbf{x}
=
\sum_{j} \int_{\hat\Omega}\hat\varphi_i(\boldsymbol{\xi})\hat\varphi_j(\boldsymbol{\xi})\,|\det J(\boldsymbol{\xi})|\,\mathrm{d}\boldsymbol{\xi}>0.
\end{equation}
\begin{proof}
Since the forward mapping $\phi:\hat\Omega\to\Omega$ must be bijective for its inverse $\phi^{-1}:\Omega\to\hat\Omega$ to exist, the determinant of $J=D\Phi$ is unconditional non-zero. The strict positivity of B-Spline basis function over their support completes the proof. \qed
\end{proof}
\end{theorem}

\textbf{Galerkin flux decomposition} \cite{KBNII}. The contribution of the residual $R(U)$ to a single DOF, say, $i$ can be decomposed into a sum of solution differences between $U_i$ and $U_j$, where $j$ extends over all DOFs for which the mass matrix satisfies $m_{ij}\ne 0$:
\begin{equation}
R_i
=
\sum_{j} \mathbf{c}_{ji}\cdot\mathbf{F}_j
=
\sum_{j\ne i} \mathbf{e}_{ij}\cdot(\mathbf{F}_j-\mathbf{F}_i)
=
\sum_{j\ne i} \mathbf{e}_{ij}\cdot\mathbf{A}_{ij}(U_j-U_i).
\label{eq:fluxdiff}
\end{equation}
Here, $\mathbf{e}_{ij}=0.5(\mathbf{c}_{ji}-\mathbf{c}_{ij})$ and $\mathbf{A}_{ij}=A(U_i,U_j)$ denotes the flux-Jacobian matrix $\mathbf{A}(U)=\partial\mathbf{F}(U)/\partial U$ evaluated for the density averaged Roe mean values \cite{Roe1981357}. 

The derivation procedure utilizes the partition of unity property of basis $\{\varphi_j\}$, which remains valid for tensor-product as well as THB-Splines \cite{GIANNELLI2012485} thus enabling adaptive refinement. 
The underlying tensor-product construction makes it possible to fully unroll the '$j\ne i$ loop' by exploiting the fact that the support of univariate B-Spline functions of order $p$ extends over a knot span of size $p+1$.

%
%
%
%
%
%
%
%


\textbf{Artificial viscosities} \cite{KBNII}. Expression \eqref{eq:fluxdiff} is augmented by  artificial viscosities
\begin{equation}
\tilde R_i:=R_i+\sum_{j\ne i}D_{ij}(U_j-U_i),
\quad
D_{ij}:=\|\mathbf{e}_{ij}\|\, R_{ij}\, |\Lambda_{ij}|\, R_{ij}^{-1}
\end{equation}
where $\Lambda_{ij}$ and $R_{ij}$ are matrices of eigenvalues and eigenvectors of $\mathbf{\tilde A}_{ij}$, respectively.

\textbf{Linearized FCT algorithm} \cite{Kuzmin2009}. Row-sum mass lumping and application of stabilizing artificial viscosities lead to the positivity-preserving predictor scheme
\begin{equation}
M_L\frac{d\tilde U}{dt}
=
\tilde R(\tilde U),
\end{equation}
which is advanced in time by the SSP-RK procedure \eqref{eq:rk3_1}--\eqref{eq:rk3_3} starting from the end-of-step solution $U^n$ of the previous time step. The solution $\tilde U^{n+1}$ and the approximation $\dot U^{n+1}\approx M_L^{-1}\tilde R(\tilde U^{n+1})$ are used to linearize the antidiffusive fluxes
\begin{equation}
F_{ij}
=
m_{ij}
\left(
\dot U_i^{n+1}-\dot U_j^{n+1}
\right)
+
\tilde D_{ij}^{n+1}
\left(
\tilde U_i^{n+1}-\tilde U_j^{n+1}
\right),
\quad
F_{ji}=-F_{ij},
\end{equation}
where $\tilde D_{ij}^{n+1}=D_{ij}(\tilde U_i^{n+1},\tilde U_j^{n+1})$ stands for the evaluation of the viscosity operator at the density-averaged Roe mean values based on the predicted solution.

The philosophy of FCT schemes is to multiply $F_{ij}$ by a symmetric correction factor $0\le \alpha_{ij}\le 1$ to obtain the amount of constrained antidiffusive correction
\begin{equation}
\bar F_i=\sum_{j\ne i}\alpha_{ij}F_{ij},
\label{eq:limit}
\end{equation}
which can be safely added to the coefficients of the positivity-preserving predictor
\begin{equation}
U^{n+1}_i=\tilde U_i^{n+1} + \frac{\Delta t}{m_i} \bar F_i
\end{equation}
without generating spurious oscillations in the updated end-of-step solution.
The flux limiting procedure developed in \cite{Kuzmin2010} calculates individual correction factors $\alpha_{ij}^u$ for user-definable scalar control variables $u(U)$, e.g., density and pressure since $\rho >0$ and $p>0$ implies $\rho E>0$ following directly from the equation of state \eqref{eq:eos}. A safe choice for the final correction factor in \eqref{eq:limit} is to set $\alpha_{ij}=\min\{\alpha_{ij}^\rho,\alpha_{ij}^p\}$.

It should be noted that the original limiting procedure has been designed for low-order \emph{nodal} finite elements, where it ensures that the nodal values of the end-of-step solution are bounded from below and above by the local minimal/maximal nodal values of the positivity-preserving predictor. That is, the following holds for all $i$:
\begin{equation}
\tilde u_i^{\min}
:=
\min_{j\ne i} \tilde u_j
\le
u_i^{n+1}
\le
\max_{j\ne i} \tilde u_j
=:
\tilde u_i^{\max}.
\end{equation}
The piece-wise linearity of $P_1$ finite elements ensures, that the corrected solution does not exceed the imposed bounds inside the elements. This is not the case for higher-order Lagrange FEM, but can be easily shown for B-Spline based FEM:
\begin{theorem}
If the flux limiter \cite{Kuzmin2010} is applied to the weights of the B-Spline expansion
\begin{equation}
\tilde u(\mathbf{x})=\sum_j \tilde u_j \varphi_j(\mathbf{x})
\end{equation}
then the flux corrected end-of-step solution (of the control variable $u$)
\begin{equation}
u^{n+1}(\mathbf{x})=\sum_j u_j^{n+1} \varphi_j(\mathbf{x})
\end{equation}
remains bounded from below and above by the bounding functions
\begin{equation}
\tilde u^{\max\atop\min}
(\mathbf{x})
=
\sum_j \tilde u_j^{\max\atop\min} \varphi_j(\mathbf{x}).
\end{equation}
\begin{proof}
Assume that for some $\mathbf{x}^*\in\Omega$ we have that $u^{n+1}(\mathbf{x}^*) > \tilde u^{\max}(\mathbf{x}^*)$. Then
\begin{equation}
0>
\tilde u^{\max}(\mathbf{x}^*)-u^{n+1}(\mathbf{x}^*)
=
\sum_j
\underbrace{[\tilde u_j^{\max}-u_j^{n+1}]}_{\ge 0}
\underbrace{\varphi(\mathbf{x}^*)}_{\ge 0}
\ge 0,
\end{equation}
which is a contradiction. A similar argument holds for the lower bound. \qed
\end{proof}
\end{theorem}

\section{Numerical results}

The proposed IGA-AFC approach has been applied to Sod's two dimensional shock tube problem \cite{Sod1978}, which has been solved on the unit square domain and the VKI U-bend test case geometry proposed in \cite{Verstraete2016}. In both cases, tensor-product bi-quadratic B-Splines combined with the third-order SSP-Runge Kutta time stepping scheme have been employed. Density and pressure have been adopted as control variables and the final correction factors have been computed as their minimum.

\begin{figure}
\begin{minipage}{0.44\textwidth}
\includegraphics[width=\textwidth]{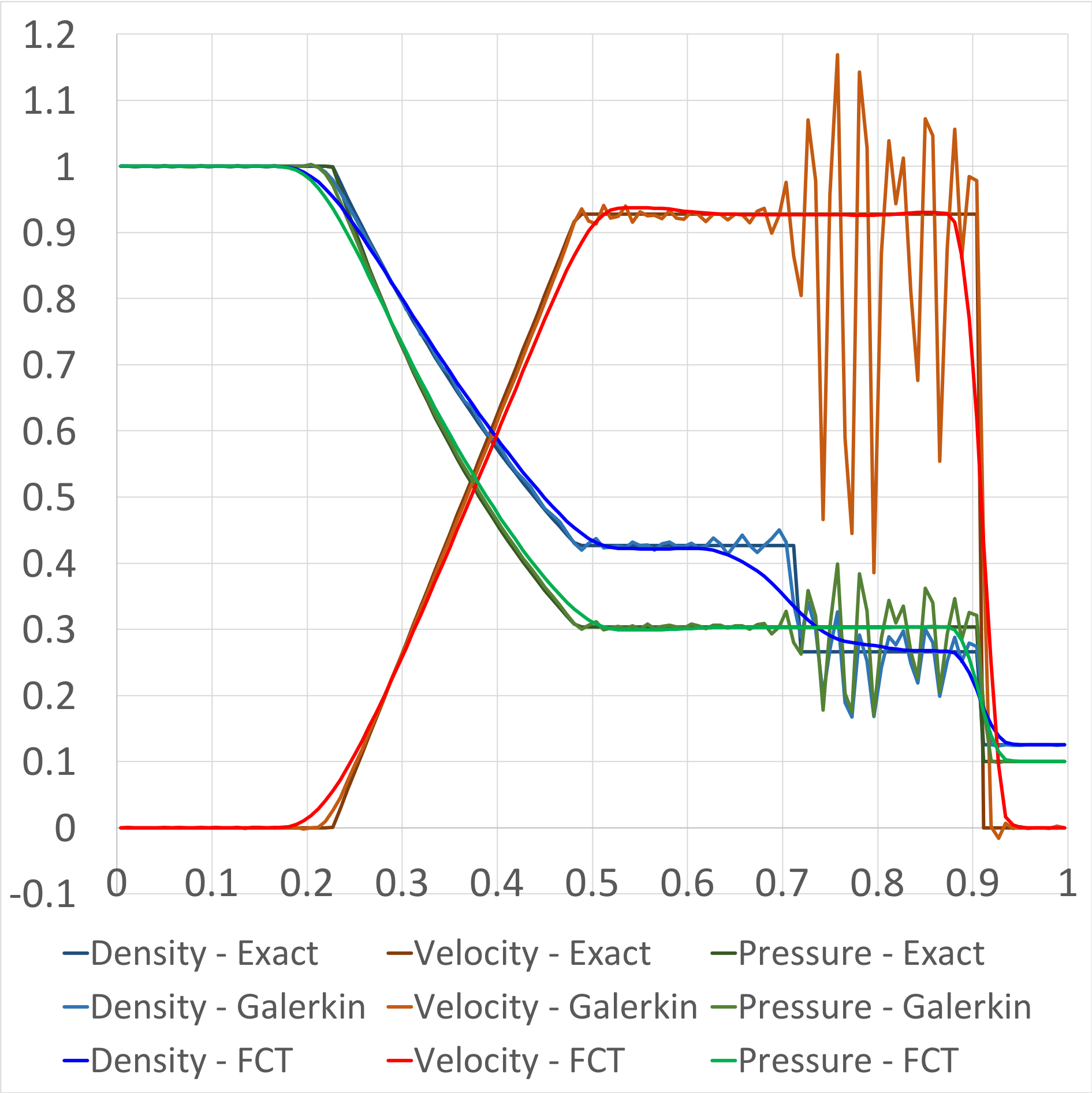}
\end{minipage}\hfill
\begin{minipage}{0.48\textwidth}
\includegraphics[width=\textwidth, trim=200 0 200 0,clip]{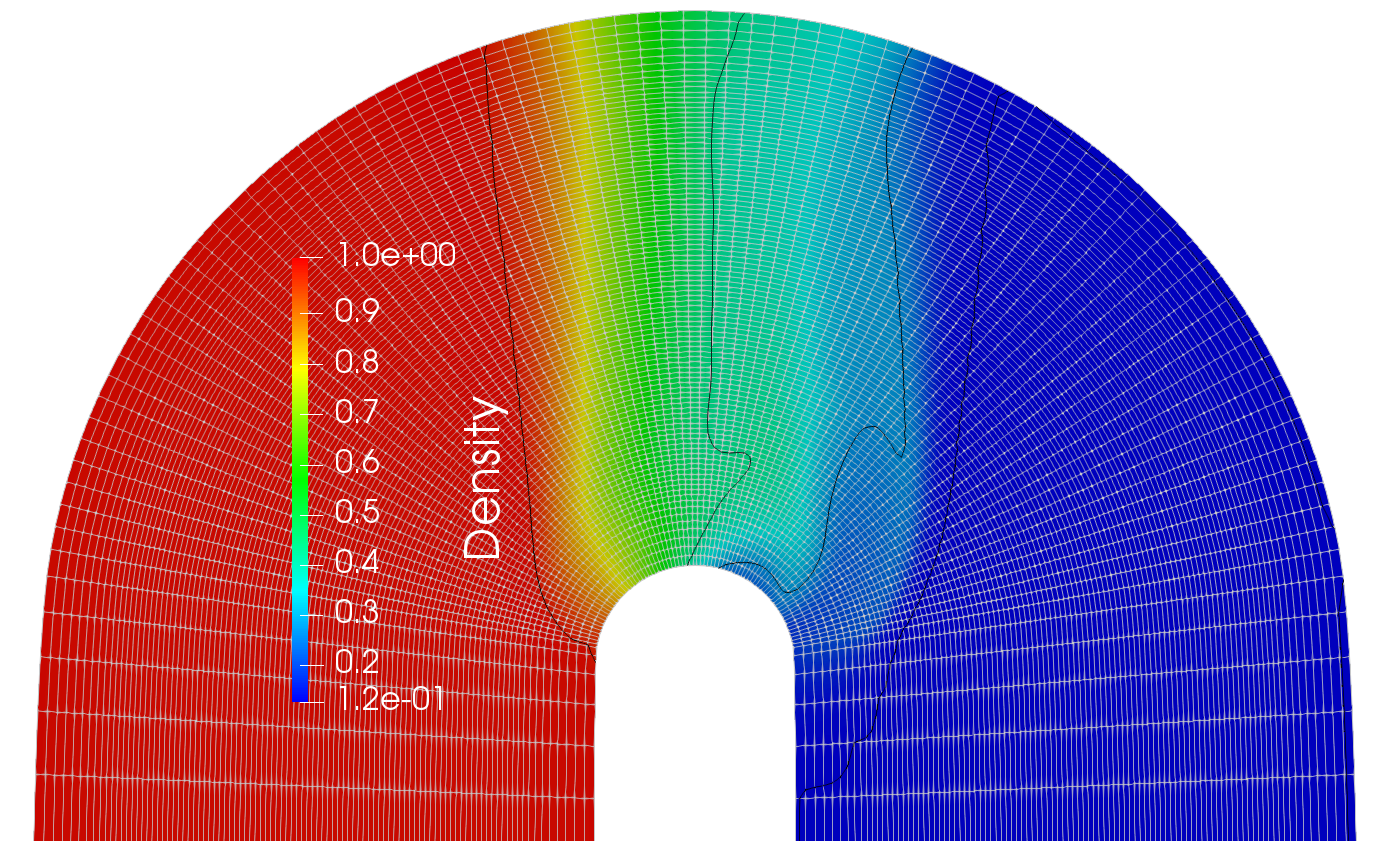}
\end{minipage}
\caption{Numerical solution to Sod's shock tube problem at $T=0.231$ computed on the unit square (left) and the VKI U-bend test case geometry (right) with bi-quadratic B-Spline basis functions.}
\label{fig:sod_shocktube}
\end{figure}

\textbf{Benchmark I}. Figure~\ref{fig:sod_shocktube} (left) shows the numerical solution sampled along the line $y=0.5$. The solution was computed on the unit square, which was discretized by $66\times 66$ equidistantly distributed B-Spline basis functions and marched forward in time with time step size $\Delta t=0.0005$ until the final time $T=0.231$ was reached. This small time-step size was needed to prevent the highly oscillatory Galerkin scheme from breaking down completely. The density and pressure profiles that were computed by IGA-AFC stay within the physical bounds and show a crisp resolution of the shock wave and the expansion fan. However, the contact discontinuity is smeared over several layers, which needs to be improved in forthcoming research.

\textbf{Benchmark II}. Figure~\ref{fig:sod_shocktube} (right) shows part of the density profile that is obtained from simulating Sod's shock tube problem on the VKI U-bend geometry \cite{Verstraete2016} with all other settings remaining unchanged except for the time step size $\Delta t=0.001$. The black contour lines indicate the location of the three characteristic wave types, i.e., the rarefaction wave, the contact discontinuity and the shock wave from left to right. As in the rectangular case, the IGA-AFC scheme succeeds in preserving the positivity of the control variables ($\rho$ and $p$, and consequently $\rho E$) , thereby demonstrating its practical applicability on non-equidistant curved 'meshes'. For illustration purposes the computational mesh has been approximated and illustrated by white lines to give an impression of the mesh width, which varies locally. 

It should be noted that the $C^1$ continuous parameterization of the curved boundary by a quadratic B-Spline function yields a unique definition of the outward unit normal vector $\mathbf{n}(\mathbf{x})=\mathbf{n}(\phi(\boldsymbol{\xi}))$ in every point on the boundary $\Gamma$. Moreover, $\mathbf{n}$ is a continuous function of the boundary parameter values $\boldsymbol{\xi}$. As a consequence, the approximate solution is free of 'numerical artifacts', which often occur for polygonal boundary representations in $P_1/Q_1$ finite elements, where the $C^0$ 'kinks' between two consecutive boundary segments serve is microscopic compression or expansion corners and, moreover, give rise to undetermined normal vectors in the nodes. In practice, this problem is often overcome by averaging the local normal vectors of the adjacent boundary segments but, still, the so-defined normal vector is not a continuous function of the boundary parameterization, which is the case for IGA.

\section{Conclusions}

In this work, we have extended the high-resolution isogeometric scheme presented in \cite{igaAFC1_fef2017} to systems of conservation laws, namely, to the compressible Euler equations. The main contribution is the positivity proof of the linearized FCT algorithm for B-Spline based discretizations, which provides the theoretical justification of our IGA-AFC approach. Future work will focus on reducing the diffusivity of the flux limiter and using THB-Splines \cite{GIANNELLI2012485} to perform adaptive refinement

\begin{acknowledgement}
This work has been supported by the European Unions Horizon 2020 research and innovation programme under grant agreement No. 678727.
\end{acknowledgement}

\bibliographystyle{unsrt}
\bibliography{bib}
\end{document}